\documentclass[a4paper,11pt]{article}
\usepackage[top=30truemm,bottom=30truemm,left=30truemm,right=30truemm]{geometry}
\usepackage{amsmath}
\usepackage{amssymb}
\usepackage{float}
\usepackage{color}
\usepackage{graphicx}
\usepackage{graphics}
\usepackage{graphicx}
\usepackage{tcolorbox}
\usepackage{latexsym}
\usepackage{ascmac}
\usepackage{framed}
\usepackage{fancybox}
\usepackage{theorem}
\usepackage{stmaryrd}
\usepackage{here}
\usepackage{mathtools}
\usepackage{bm}
\usepackage{appendix}
\usepackage{url}

\mathtoolsset{showonlyrefs=true} 

\theoremstyle{break}
\def\qed{\hfill$\Box$}

\newtheorem{defn}{Definition}
\newtheorem{lem}{Lemma}[section]
\newtheorem{prop}{Proposition}
\newtheorem{thm}{Theorem}
\newtheorem{cor}{Corollary}

\newtheorem{ex}{Example}
\newtheorem{rem}{Remark}

\makeatletter
\@addtoreset{equation}{section}

\makeatother

\title{{\bf{
Criteria for memory-induced blow-up in distributed delay differential equations
}}
}

\author{
Yu Ichida
\thanks{Department of Mathematics, School of Science and Technology, Meiji University, 1-1-1 Higashimita, Tama-ku, Kawasaki 214-8571, Japan, {\tt ichidayu@meiji.ac.jp} }
}

\begin{document}
\maketitle

\begin{abstract}
This paper investigates the conditions under which solutions to distributed delay differential equations describing memory accumulation and forgetting blow up in finite time. 
It provides criteria for memory-induced blow-up, which means finite-time blow-up of the delayed equation while the corresponding non-delayed equation, obtained in a certain limit, has global solutions. 
In the author's previous work, memory-induced blow-up was shown for specific nonlinearities by means of Poincar\'e-type compactification, but a general theory was not established. 
Here, the two-dimensional ODE system derived by the linear chain trick is rewritten as an orbit equation, and a comparison theorem and an Osgood-type integral condition are applied to obtain a blow-up criterion based on the existence of a subsolution. 
This characterizes the nonlinearities for which memory-induced blow-up occurs. 
This paper also gives a necessary condition in a special case, a global existence condition based on a supersolution, and a systematic method for constructing subsolutions for separable nonlinearities. 
These results recover the author's previous results and are illustrated through several examples.
\end{abstract}


\begin{center}
\end{center}

\section{Introduction}
\label{sec:DDDEg-int}
This paper is concerned with finite-time blow-up of solutions to the following distributed delay differential equation:
\begin{equation}
\dot{u}(t) = f\left( u, \int_{-\infty}^{t} \alpha e^{-\beta(t-s)} g(u(s)) \, ds \right).
\label{eq:DDDEg-int1}
\end{equation}
Here $u=u(t)$ is the unknown function, $t$ denotes the present time, $s$ denotes the past ($-\infty<s\le t$), $\alpha$ and $\beta$ are positive constants, and $\dot{u}=du/dt$. 
The parameter $\alpha$ represents the strength of history or memory, while $\beta$ corresponds to the rate of forgetting. 
If $\beta$ is large, memory decays quickly and the effect of the delay is weak; if $\beta$ is small, the influence of memory, that is, of the time delay, is strong. 
Let $\mathbb{R}_+$ denote the set of positive real numbers. 
We impose the initial function as follows:
\begin{equation}
u(s)=\phi(s), \quad \phi \in C_{b}:=C_{b}((-\infty, 0]; \mathbb{R}_{+}),
\quad \phi(s)>0,
\quad -\infty < s \le 0.
\label{eq:DDDEg-int2}
\end{equation}
Here $\phi\in C_b((-\infty,0];\mathbb{R}_+)$ means that $\phi$ is a positive continuous function from $(-\infty,0]$ to $\mathbb{R}_+$ satisfying the boundedness condition $\sup_{s\le 0}|\phi(s)|<+\infty$.

We further define $v$ by
\begin{equation}
v(t)=\int_{-\infty}^{t} \alpha e^{-\beta(t-s)} g(u(s)) \, ds.
\label{eq:DDDEg-int3}
\end{equation}
By the choice of the initial function \eqref{eq:DDDEg-int2} and $\beta>0$, the total amount of history and memory $v(0)$ at the initial time $t=0$ is finite. 
Moreover, if $g(\phi(0))>0$, then $v(0)>0$.

One of the central questions in the qualitative theory of delay differential equations is how time delays change the behavior of solutions. 
As noted in \cite{DDDEb}, time delays can destabilize solutions, induce oscillations, or, conversely, stabilize solutions (see, e.g., Smith \cite{Smith} and the references therein).

We now review representative studies on finite-time singularities in delay differential equations, following the survey given in \cite{DDDEb}. 
Ezzinbi and Jazar \cite{EJ2006} studied scalar differential equations with a constant delay and clarified, by comparison with ordinary differential equations, how additive and multiplicative delay terms induce or suppress finite-time blow-up. 
To the best of the author's knowledge, the first result in which a time delay essentially induces blow-up is due to Eremin, Ishiwata, Ishiwata, and Nakata \cite{EIIN2021}. 
They took the normal form of the Hopf bifurcation as the non-delayed system, introduced a constant delay, and demonstrated delay-induced blow-up. 
Subsequently, Yagasaki \cite{Yagasaki2021} used the averaging method, starting from the normal form of the subcritical Hopf bifurcation, to show finite-time blow-up of solutions under suitable conditions. 
The effect of a constant delay in a scalar equation was also investigated by Ishiwata and Nakata \cite{IN2021}. 
However, to the best of the author's knowledge, there are almost no results on blow-up of solutions caused by distributed delays, which describe the accumulation of past history and memory.

Finite-time blow-up of solutions to Volterra-type integro-differential equations and integral equations has also been studied in recent years \cite{appl,bruya}.
It has been shown that the singularity of the kernel and its behavior near zero govern the blow-up behavior. 
On the other hand, mathematical models in the natural and life sciences that involve the accumulation and forgetting of history or memory are frequently described not by mere time-integral equations but by distributed delay differential equations, in which differentiation and integration are combined. 
The mechanism by which delays or memory essentially induce new blow-up has not been sufficiently discussed, and the pioneering study on finite-time blow-up of solutions to differential equations with distributed delays is the author's work (\cite{DDDEb}).

In \cite{DDDEb}, the author showed that solutions $u=u(t)$ of \eqref{eq:DDDEg-int1} blow up in finite time for the following choices of $f$ and $g$:
\begin{itemize}
  \item[(F1)] $f(u,v)=v^{p}$ ($p>1$) and $g(u)=u$,
  \item[(F2)] $f(u,v)=(1-v)^{-1}$ and $g(u)=u$,
  \item[(F3)] $f(u,v)=(1-v)^{-1}-\gamma u$ and $g(u)=u-a$, where $\gamma$ and $a$ are positive constants.
\end{itemize}
This paper aims to build a general theory applicable to the cases (F1) and (F2). 
The case (F3) exhibits a threshold phenomenon for blow-up, which was shown in \cite{DDDEb}; it is excluded from the present paper in order to keep the theoretical framework clear.

The idea of the proof in \cite{DDDEb} is as follows. 
Defining $v=v(t)$ by \eqref{eq:DDDEg-int3} for the delay differential equation \eqref{eq:DDDEg-int1} with the initial condition \eqref{eq:DDDEg-int2}, we obtain
\begin{equation}
\dot{v}= -\beta v(t)+\alpha g(u(t)). 
\label{eq:DDDEg-int4}
\end{equation}
Hence \eqref{eq:DDDEg-int1} reduces to the following two-dimensional system of ordinary differential equations:
\begin{equation}
\begin{cases}
\dot{u}=f(u, v), \\
\dot{v}=  \alpha g(u(t)) -\beta v(t).
\end{cases}
\label{eq:DDDEg-int5}
\end{equation}
This idea of reducing a delay differential equation to a system of ordinary differential equations is called the linear chain trick (see, e.g., \cite{MacDonald1978, MacDonald1989, Ruan2006, Smith} and the references therein). 
In \cite{DDDEb}, the dynamics of \eqref{eq:DDDEg-int5} on the phase space including infinity was analyzed by combining the Poincar\'e-type compactification, the blow-up technique for singular points, and the center manifold theorem. 
This analysis showed that solutions of \eqref{eq:DDDEg-int1} blow up in finite time, and the blow-up rate was derived.

Furthermore, in \cite{DDDEb}, memory-induced blow-up, that is, blow-up essentially induced by memory, was demonstrated explicitly by comparing the delay equation with the non-delayed equations obtained in the limit $\beta\to+\infty$ (see Definition \ref{def:DDDEg-mr1} below for the definition of memory-induced blow-up). 
If $g(u(\cdot))$ is sufficiently well-behaved, then by Lebesgue's dominated convergence theorem, for fixed $c:=\alpha/\beta$,
\begin{equation}
\lim_{\beta\to +\infty} v(t) = cg(u(t)).
\label{eq:DDDEg-int6}
\end{equation}
For fixed $\alpha$, $v(t)\to 0$ as $\beta\to+\infty$. 
Therefore, the (non-delayed) equations obtained as $\beta\to+\infty$ are the following:
\begin{enumerate}
  \item[(i)] for fixed $c$ (in particular, $c=1$, i.e., $\alpha=\beta$): $\dot{u}=f(u,c\,g(u))$;
  \item[(ii)] for fixed $\alpha$ ($c\to0$): $\dot{u}=\hat{f}(u):=f(u,0)$.
\end{enumerate}
Accordingly, for $c=1$, the non-delayed equations corresponding to (F1) and (F2) are
\begin{equation}
{\rm{(F1)}}:\ \dot{u}(t)= u(t)^{p}, \quad
{\rm{(F2)}}:\ \dot{u}(t)= \dfrac{1}{1-u(t)},
\label{eq:DDDEg-int7}
\end{equation}
which are the typical examples of blow-up and quenching, well-known finite-time singularities of ordinary differential equations. 
In contrast, for $c\to0$, the non-delayed equations corresponding to (F1) and (F2) are
\begin{equation}
{\rm{(F1)}}:\ \dot{u}(t)= 0, \quad
{\rm{(F2)}}:\ \dot{u}(t)= 1.
\label{eq:DDDEg-int8}
\end{equation}
In both cases the solutions of the non-delayed equations exist globally. 
Hence, the distributed delay with memory effect essentially contributes to the finite-time blow-up of solutions; that is, memory-induced blow-up occurs.

Motivated by these results of \cite{DDDEb}, the following question arises naturally: for which $f$ and $g$ do solutions of the distributed delay differential equation \eqref{eq:DDDEg-int1} blow up in finite time, and for which do they exhibit memory-induced blow-up? 
Since the method of \cite{DDDEb} cannot be applied to all $f$ and $g$, it is not expected to lead to a general theory. 
It therefore suffices to look for conditions under which $u=u(t)$ in the system \eqref{eq:DDDEg-int5}, obtained by the linear chain trick, blows up in finite time. 
For this purpose, the Osgood-type condition (Proposition \ref{prop:DDDEg-pre4}), an integral convergence test known as a criterion for finite-time blow-up of first-order ordinary differential equations, which is reviewed in Section \ref{sec:DDDEg-pre}.

However, in the first equation of \eqref{eq:DDDEg-int5}, $f$ also depends on $v(t)$, so the Osgood-type condition cannot be applied in a straightforward manner. 
In other words, the key point of the analysis is to estimate and control appropriately the dynamical variation caused by the delay term (memory effect). 
In addition, this paper focuses on $v_{\max}\in(0,\infty]$, the right endpoint of the domain of $f$ with respect to $v$ (the location of the singularity). 
When $f(u,v)$ has a singularity such that $f(u,v)\to\infty$ as $v\to v_{\max}^{-}$, as in the case of $(1-v)^{-1}$ in (F2), finite-time singularities may occur not only through the divergence of the unknown function $u(t)$ itself (blow-up), but also through the divergence of $f$ caused by the memory $v(t)$ reaching the limit value $v_{\max}$. 
Since the mechanism triggering blow-up is essentially different depending on whether $v_{\max}$ is finite or infinite, it is indispensable to distinguish and analyze these cases in a unified and explicit manner.

To address these issues, this paper adopts the approach of reducing the two-dimensional system \eqref{eq:DDDEg-int5}, obtained by the linear chain trick, to an orbit equation with $u$ as the independent variable. 
Then, we use a comparison theorem for ordinary differential equations (Proposition \ref{prop:DDDEg-pre3}) to construct dynamic lower bounds, subsolutions, and supersolutions.
This paper derives a sufficient condition for finite-time blow-up (Theorem \ref{th:DDDEg-mr1} and Corollary \ref{cor:DDDEg-mr1}), a necessary condition (Theorem \ref{th:DDDEg-mr2}), and a global existence condition (Theorem \ref{th:DDDEg-mr3}) for $f$ and $g$ satisfying assumptions (I) and (II) in Section \ref{sec:DDDEg-mr}. 
These results describe explicit criteria for the mechanism by which distributed delays with memory effect essentially contribute to finite-time blow-up of solutions, within the framework of memory-induced blow-up (see Definition \ref{def:DDDEg-mr1}).

In Section \ref{sec:DDDEg-ex}, as applications of the main results, finite-time blow-up, memory-induced blow-up, and global existence are determined for concrete pairs of $f$ and $g$. 
First, a procedure for systematically constructing the subsolution $\Phi$ in assumption (III-2) when $f$ is of separable form is given (Proposition \ref{prop:DDDEg-ex1}), and the results of the previous work \cite{DDDEb} for (F1) and (F2) are reproduced. 
Furthermore, the following examples are presented: an example that cannot be handled immediately by the method of \cite{DDDEb} but for which memory-induced blow-up is established by the criteria of this paper; examples in which the sufficient condition (Theorem \ref{th:DDDEg-mr1}) and the necessary condition (Theorem \ref{th:DDDEg-mr2}) coincide, yielding a necessary and sufficient condition for finite-time blow-up; and an example in which solutions blow up in finite time but the blow-up is not memory-induced. 
These show that the main results of this paper allow finite-time blow-up and global existence of solutions to equations with memory effect to be determined within a common framework based on subsolutions and supersolutions. 
Moreover, they show that condition (IV) distinguishes the case where blow-up is caused by memory (memory-induced) from the case where blow-up already arises from the non-delayed part (Example \ref{ex:DDDEg-ex6}).

The rest of this paper is organized as follows.
The next section briefly reviews the Osgood-type condition for first-order ordinary differential equations and the comparison theorem for solutions of ordinary differential equations. 
Section \ref{sec:DDDEg-mr} states the main results on finite-time blow-up and memory-induced blow-up of the solution $u=u(t)$ of \eqref{eq:DDDEg-int1}. 
Section \ref{sec:DDDEg-pro} gives rigorous proofs of the main results, and Section \ref{sec:DDDEg-ex} presents several applications of the main results of Section \ref{sec:DDDEg-mr}. 
Section \ref{sec:DDDEg-dis} gives a discussion on the comparison with non-delayed equations and some remarks.

\section{Preliminaries}
\label{sec:DDDEg-pre}
Let $F$ be a locally Lipschitz continuous function on $(0,+\infty)$, and consider the initial value problem for the ordinary differential equation:
\begin{equation}
\dot{u}=F(u(t)), \quad
u(0)=u_{0}>0.
\label{eq:DDDEg-pre1}
\end{equation}
By the local existence and uniqueness theorem and the continuation theorem for ordinary differential equations (see, e.g., \cite{ber, hara, teschl}), for every initial value $u_0>0$ there exists a unique maximal solution $u=u(t)$ defined on a maximal interval of existence $[0,T)$ with $T\in(0,\infty]$.

We say that the solution $u=u(t)$ of \eqref{eq:DDDEg-pre1} blows up in finite time if there exists $T\in(0,\infty)$ such that
\begin{equation}
\lim_{t \to T-0} u(t) = \infty.
\label{eq:DDDEg-pre2}
\end{equation}
If $T=\infty$, the solution is said to exist globally. 
As is well known, the following Osgood-type condition gives a necessary and sufficient condition for finite-time blow-up of the solution of \eqref{eq:DDDEg-pre1}.

\begin{prop}
\label{prop:DDDEg-pre4}
Let $u_0>0$, and suppose that $F$ is locally Lipschitz continuous and $F>0$ on $[u_0,\infty)$. 
Let $U$ be the maximal solution of the initial value problem \eqref{eq:DDDEg-pre1} and let $T\in(0,\infty]$ be its maximal existence time. 
Then the following holds:
\begin{equation}
 T=\int_{u_0}^{\infty} \dfrac{d\xi}{F(\xi)} \in (0,\infty].
 \label{eq:DDDEg-pre18}
\end{equation}
In particular,
\begin{enumerate}
\item[(i)]
the solution $U$ exists globally on $[0,\infty)$ if and only if $\int_{u_0}^{\infty}d\xi/F(\xi)=\infty$;
\item[(ii)] 
the solution $U$ blows up in finite time if and only if $\int_{u_0}^{\infty}d\xi/F(\xi)<\infty$. In this case, $\lim_{t\to T-0}U(t)=\infty$, that is, \eqref{eq:DDDEg-pre2} holds.
\end{enumerate}
\end{prop}

Next, to prove the main results of this paper, we introduce the following comparison theorem. 
In Section \ref{sec:DDDEg-pro}, we estimate the solution $V$ of the orbit equation with $u$ as the independent variable,
\[
V'(u)=H(u, V(u)), \quad
H(u, w):= \dfrac{\alpha g(u)-\beta w}{f(u,w)}
\]
from below and above by a subsolution $\Phi$ and a supersolution $\Psi$, respectively. 
For a general $H$, the following comparison theorem holds.

\begin{prop}
\label{prop:DDDEg-pre3}
Let $u_1<\omega\le\infty$ and $v_{\max}\in(0,\infty]$, and set $I:=[u_1,\omega)$ and $J:=(0,v_{\max})$. 
Suppose that $H:I\times J\to\mathbb{R}$ is continuous and locally Lipschitz continuous in the second variable; that is, for every compact set $K\subset I\times J$ there exists a constant $L_K\ge0$ such that
\begin{equation}
  |H(u,w)-H(u,\tilde w)|\le L_K\,|w-\tilde w|
  \quad \bigl((u,w),(u,\tilde w)\in K\bigr).
  \label{eq:DDDEg-pre15}
\end{equation}
If $\Phi,\Psi\in C^1(I;J)$ satisfy
\begin{equation}
  \Phi'(u)\le H\bigl(u,\Phi(u)\bigr),\quad
  \Psi'(u)\ge H\bigl(u,\Psi(u)\bigr)\quad(u\in I),\quad
  \Phi(u_1)\le\Psi(u_1),
   \label{eq:DDDEg-pre16}
\end{equation}
then $\Phi(u)\le\Psi(u)$ on $I$.
\end{prop}

\begin{rem}
\label{rem:DDDEg-pre1}
In the proof of Lemma \ref{lem:DDDEg-pro3}, Proposition \ref{prop:DDDEg-pre3} is applied with $(\Phi,\Psi)=(\Phi,V)$, regarding the solution $V$ as a supersolution (equality holds). 
On the other hand, in the proof of Theorem \ref{th:DDDEg-mr3}, Proposition \ref{prop:DDDEg-pre3} is applied with $(\Phi,\Psi)=(V,\Psi)$, regarding the solution $V$ as a subsolution (equality holds).
\end{rem}

\section{Main results}
\label{sec:DDDEg-mr}
Before stating the main results, we define memory-induced blow-up for the distributed delay differential equation \eqref{eq:DDDEg-int1}. 
Hereafter, $u_0=\phi(0)$ by \eqref{eq:DDDEg-int2}.

\begin{defn}
\label{def:DDDEg-mr1}
Assume that $f$ extends continuously to $v\to0^+$ and that $f(\cdot,0)$ is locally Lipschitz continuous. 
We say that a solution $u$ of the distributed delay differential equation \eqref{eq:DDDEg-int1}, \eqref{eq:DDDEg-int2} exhibits memory-induced blow-up if $u$ satisfies \eqref{eq:DDDEg-pre2} at some finite time $T<\infty$, and the solution $U=U(t)$ of the non-delayed equation
\[
\dot{U}=\hat f_0(U)=f(U,0), \quad U(0)=u_{0}
\]
obtained in the limit $c:=\alpha/\beta\to0$, exists globally on $[0,\infty)$.
\end{defn}

Let $v_{\max}\in(0,\infty]$, and let $\Omega_+=(0,\infty)\times(0,v_{\max})$ be the domain of $(u,v)$. 
For the initial function $\phi\in C_b$ given in \eqref{eq:DDDEg-int2}, set
\[
v(0):=\alpha\int_{-\infty}^0e^{\beta s}g(\phi(s))\,ds
\]
and assume $v(0)\in(0,v_{\max})$. 
Furthermore, define $v_{\min}$ by
\begin{equation}
v_{\min} := \min\left\{ v(0), \, \frac{\alpha}{\beta} g(u_0) \right\} >0.
\label{eq:DDDEg-mr1}
\end{equation}
In this paper, we impose the following assumptions on $f$ and $g$.
\begin{enumerate}
\item[(I)] 
$f$ is of class $C^1$ on $\Omega_+$ and satisfies $f>0$ and $f_v\ge0$ for $v>0$. 
Moreover, $f$ extends continuously to $v=0$, and $f(u,0)$, understood as this extended value, is locally Lipschitz continuous.
\item[(II)] $g$ is of class $C^1$ on $[0,+\infty)$ and satisfies $g\ge0$, $g'(u)\ge0$, $g\not\equiv0$, and $g(u_0)>0$.
\item[(III-1)]
Let $v(0)\in(0,v_{\max})$ and let $v_{\min}$ be as in \eqref{eq:DDDEg-mr1}. 
Then the following holds
\begin{equation}
\int_{u_{0}}^{+\infty}\dfrac{du}{ f(u, v_{\min})} <\infty.
\label{eq:DDDEg-mr6}
\end{equation}
\item[(III-2)] 
Let $v(0)\in(0,v_{\max})$ and let $v_{\min}$ be as in \eqref{eq:DDDEg-mr1}. 
Then there exist $u_1\ge u_0$ and a function $\Phi\in C^1([u_1,\infty);(0,v_{\max}))$ such that the followings hold:
\begin{equation}
\Phi'(u) \le \frac{\alpha g(u)-\beta\Phi(u)}{f(u,\Phi(u))} \quad (u\ge u_1),\quad \Phi(u_1)\le v_{\min}
\label{eq:DDDEg-mr2}
\end{equation}
and 
\begin{equation}
\int_{u_1}^{+\infty} \frac{du}{f(u,\Phi(u))}<\infty.
\label{eq:DDDEg-mr3}
\end{equation}
\item[(IV)] 
The solution $U$ of $\dot U=f(U,0)$, $U(0)=u_0$, exists globally on $[0,\infty)$. 
When $f(\cdot,0)>0$, by Proposition \ref{prop:DDDEg-pre4} this is equivalent to
\begin{equation}
\int_{u_0}^{\infty} \frac{du}{f(u,0)} = \infty.
\label{eq:DDDEg-mr4}
\end{equation}
\item[(V)]
$v_{\max}<\infty$, and there exist $\kappa>0$, $q\ge1$, and $\delta_0\in(0,v_{\max})$ such that for all $u\ge u_0$ and $v\in(v_{\max}-\delta_0,v_{\max})$,
\begin{equation}
f(u,v) \ge \kappa(v_{\max}-v)^{-q}
\label{eq:DDDEg-mr5}
\end{equation}
holds.
In particular, $f\to\infty$ as $v\to v_{\max}^-$.
\end{enumerate}

Then, we now state the main results in this paper.

\begin{thm}[Sufficient condition]
\label{th:DDDEg-mr1}
Let $\alpha$ and $\beta$ be positive constants, and assume (I), (II), and (III-2). 
If $v_{\max}\in(0,\infty)$, assume also (V). 
Then the solution $u=u(t)$ of the distributed delay differential equation \eqref{eq:DDDEg-int1}, \eqref{eq:DDDEg-int2} has a maximal existence time $T_{\max}\in(0,+\infty)$ and satisfies
\begin{equation}
\lim_{t\to T_{\max}^{-}} u(t) =+\infty.
\label{eq:DDDEg-mr11}
\end{equation}
Moreover, the following holds:
\begin{equation}
 T_{\max} \le\frac{u_1-u_0}{c_1}+\int_{u_1}^\infty\frac{du}{f(u,\Phi(u))},\quad
 c_1:=\min_{u_0\le u\le u_1}f(u,v_{\min})>0.
\label{eq:DDDEg-mr8}
\end{equation}
If, in addition, (IV) holds, then this blow-up is memory-induced blow-up.
\end{thm}

Regarding assumption (III-2), taking $\Phi(u)\equiv v_{\min}$ yields the following corollary.

\begin{cor}
\label{cor:DDDEg-mr1}
Let $\alpha$ and $\beta$ be positive constants, and assume (I), (II), and (III-1). 
If $v_{\max}\in(0,\infty)$, assume also (V). 
Then the conclusion of Theorem \ref{th:DDDEg-mr1} holds.
\end{cor}

Theorem \ref{th:DDDEg-mr1} and Corollary \ref{cor:DDDEg-mr1} give sufficient conditions for the solution $u=u(t)$ of the distributed delay differential equation \eqref{eq:DDDEg-int1}, \eqref{eq:DDDEg-int2} to blow up in finite time and to exhibit memory-induced blow-up. 
It is natural to ask how sharp these conditions are. 
The following theorem gives a necessary condition when $v_{\max}=\infty$.

\begin{thm}[Necessary condition]
\label{th:DDDEg-mr2}
Let $v_{\max}=\infty$ and assume (I) and (II). 
If the solution of \eqref{eq:DDDEg-int1}, \eqref{eq:DDDEg-int2} blows up at a finite time $T_{\max}<\infty$, then
\begin{equation}
 \int_{u_0}^\infty\frac{du}{f\bigl(u,\ \max\{v(0),\tfrac{\alpha}{\beta}g(u)\}\bigr)}\le T_{\max}<\infty.
 \label{eq:DDDEg-mr7}
\end{equation}
Equivalently, if the integral on the left-hand side is $+\infty$, then the solution exists globally.
\end{thm}

\begin{thm}
\label{th:DDDEg-mr3}
Let $v_{\max}\in(0,\infty]$ and assume (I) and (II). 
If $\Psi\in C^1([u_0,\infty);(0,v_{\max}))$ satisfies
\begin{equation}\label{eq:DDDEg-mr9}
 \Psi(u_0)\ge v(0),\quad
 \Psi'(u)\ge \dfrac{\alpha g(u)-\beta\Psi(u)}{f(u,\Psi(u))}\quad(u\ge u_0),\quad
 \int_{u_0}^\infty\frac{du}{f(u,\Psi(u))}=\infty
\end{equation}
then the solution of \eqref{eq:DDDEg-int1}, \eqref{eq:DDDEg-int2} exists globally (that is, $T_{\max}=\infty$).
\end{thm}

\begin{rem}
\label{rem:DDDEg-mr2}
Theorem~\ref{th:DDDEg-mr3} gives a sufficient condition for global existence in terms of a supersolution. 
Assume that 
\begin{equation}
v(0)\le\dfrac{\alpha}{\beta}g(u_0),\quad
\dfrac{\alpha}{\beta}g(u)<v_{\max} \quad (u\ge u_{0}).
\label{eq:DDDEg-mr10}
\end{equation}
Then $\Psi(u):=\alpha\beta^{-1}g(u)$ satisfies the first two conditions in \eqref{eq:DDDEg-mr9}. 
Indeed, by $g\ge g(u_0)>0$ and \eqref{eq:DDDEg-mr10}, we have $\Psi\in C^1([u_0,\infty);(0,v_{\max}))$ and $\Psi(u_0)=\alpha\beta^{-1}g(u_0)\ge v(0)$, and
\[
 \frac{\alpha g(u)-\beta\Psi(u)}{f(u,\Psi(u))}=0\le\frac{\alpha}{\beta}g'(u)=\Psi'(u).
\]
Therefore, under \eqref{eq:DDDEg-mr10}, if
\[
 \int_{u_0}^\infty\frac{du}{f\bigl(u, \alpha\beta^{-1}g(u)\bigr)}=\infty
\]
then the solution exists globally. 
When $v_{\max}=\infty$, the second condition in \eqref{eq:DDDEg-mr10} holds automatically ($g$ may be unbounded). 
Theorem \ref{th:DDDEg-mr2} can be used regardless of the size of $v(0)$, but is restricted to $v_{\max}=\infty$. 
On the other hand, Theorem \ref{th:DDDEg-mr3} can also be used when $v_{\max}<\infty$.
\end{rem}

\section{Proof of Theorems}
\label{sec:DDDEg-pro}
In this section, we prove Theorems~\ref{th:DDDEg-mr1}, \ref{th:DDDEg-mr2}, \ref{th:DDDEg-mr3} and Corollary~\ref{cor:DDDEg-mr1}. 
As mentioned in Section \ref{sec:DDDEg-int}, for the distributed delay differential equation \eqref{eq:DDDEg-int1}, setting $v=v(t)$ as in \eqref{eq:DDDEg-int3} yields, by the linear chain trick (see, e.g., \cite{MacDonald1978, MacDonald1989, Ruan2006, Smith} and the references therein), the two-dimensional system of ordinary differential equations \eqref{eq:DDDEg-int5}. We consider \eqref{eq:DDDEg-int5} on the domain
\begin{equation}
\Omega_+ = (0, \infty) \times (0, v_{\max}), \quad
v_{\max}\in (0, +\infty].
\label{eq:DDDEg-pro1}
\end{equation}
Since $f\in C^1(\Omega_+)$ and $g\in C^1([0,\infty))$ are assumed, the right-hand side is locally Lipschitz continuous. 
Hence for every initial value $(u_0,v(0))\in(0,\infty)\times(0,v_{\max})$ there exists a unique local solution $(u,v)\in C^1([0,T_{\max}))^2$ on a maximal interval of existence $[0,T_{\max})$, where $T_{\max}\in(0,\infty]$ is the maximal existence time. 
Moreover, the continuation theorem holds: if $T_{\max}<\infty$, then $(u(t),v(t))$ cannot remain in a compact set $K\subset\Omega_+$ for all $t\in[0,T_{\max})$.

The strategy of the proofs is as follows.
We apply the comparison theorem to the orbit equation with $u$ as the independent variable to show $T_{\max}<\infty$, and then show that $u(t)\to\infty$ under $T_{\max}<\infty$.

\subsection{Monotonicity and a lower bound for the memory}
\label{sub:DDDEg-pro1}

\begin{lem}
\label{lem:DDDEg-pro1}
Assume (I) and (II). 
Then the solution $(u,v)$ satisfies, on $[0,T_{\max})$, the following holds
\begin{equation}
\dot{u}>0,\quad
u(t) \ge u_0, \quad 
0<v_{\min} \le v(t) \le \max\left\{ v(0), \dfrac{\alpha}{\beta}g(u(t)) \right\},
\label{eq:DDDEg-pro3}
\end{equation}
where $v_{\min}$ is defined in \eqref{eq:DDDEg-mr1}.
\end{lem}
{\bf{Proof.}}
Since $f(u,v)>0$ for $v>0$ by assumption, we have $\dot u=f(u,v)>0$. 
Hence $u$ is strictly increasing and $u(t)\ge u_0$. 
This means that \eqref{eq:DDDEg-int5} has no equilibrium in $\Omega_+$. 
Next, multiplying the second equation $\dot v+\beta v=\alpha g(u)$ of \eqref{eq:DDDEg-int5} by $e^{\beta t}$ and integrating from $0$ to $t$, we obtain
\begin{equation}
 v(t) = e^{-\beta t}v(0) + \alpha \int_0^te^{-\beta(t-s)}g(u(s))\,ds.
 \label{eq:DDDEg-pro5}
\end{equation}
By $u(t)\ge u_0$ and the monotonicity of $g$, for $0\le s\le t$,
\[
g(u_0)\le g(u(s))\le g(u(t))
\]
holds.
Hence, setting $\delta:=e^{-\beta t}\in(0,1]$, we have
\begin{equation}
 \delta v(0) + (1-\delta)\frac{\alpha}{\beta}g(u_0) \le v(t) \le \delta v(0)+(1-\delta)\dfrac{\alpha}{\beta}g(u(t)).
 \label{eq:DDDEg-pro4}
\end{equation}
Both sides of \eqref{eq:DDDEg-pro4} are convex combinations of $v(0)$ and $\alpha\beta^{-1}g(u_0)$ (in the right-hand side, of $v(0)$ and $\alpha\beta^{-1}g(u(t))$). 
The case $\delta=1$ ($t=0$) corresponds to the initial value $v(0)$, the total past history, while $\delta\to0$ ($t\to\infty$) corresponds to the steady level $\alpha\beta^{-1}g$ of memory continually generated by the present state. 
Since a convex combination of two values lies between their minimum and maximum, we obtain
\begin{equation}
0<v_{\min} \le v(t) \le \max\left\{ v(0), \dfrac{\alpha}{\beta}g(u(t)) \right\}.
\label{eq:DDDEg-pro6}
\end{equation}
This gives \eqref{eq:DDDEg-pro3}, which completes the proof of Lemma \ref{lem:DDDEg-pro1}.
\qed
\\

By Lemma \ref{lem:DDDEg-pro1}, the limit
\begin{equation}
U^* := \lim_{t\to T_{\max}^{-}} u(t) \in (u_0,\infty]
\label{eq:DDDEg-pro9}
\end{equation}
exists. 
Moreover, since $v_{\min}\le v(0)<v_{\max}$, we have $v_{\min}\in(0,v_{\max})$.
By assumption (I), the following holds:
\begin{equation}
f(u, v(t)) \ge f(u, v_{\min}). 
\label{eq:DDDEg-pro10}
\end{equation}

\subsection{The orbit equation, a dynamic lower bound, and the subsolution condition}
\label{sub:DDDEg-pro2}
In the examples of the next section, the convergence test \eqref{eq:DDDEg-mr6} based on $f(u,v_{\min})$, which uses the constant $v_{\min}$ as a lower bound for $f(u,v)$, turns out to be insufficient. 
We therefore consider a dynamic lower bound for $f(u,v)$ based on the orbit equation.

Since $\dot u>0$, the map $u:[0,T_{\max})\to[u_0,U^*)$ is a $C^1$ bijection, and its inverse $t=t(u)$ is also of class $C^1$. 
Define $V(u):=v(t(u))$ for $u_0\le u<U^*$, and set
\begin{equation}
 H(u,V):=\frac{\alpha g(u)-\beta V}{f(u,V)} 
 \quad ((u,V) \in \Omega_{+}).
 \label{eq:DDDEg-pro11}
\end{equation}
Since $f>0$ and $f,g\in C^1$, we have $H\in C^1(\Omega_+)$. 
Then the following holds.

\begin{lem}
\label{lem:DDDEg-pro2}
$V$ is of class $C^1$ on $[u_0,U^*)$ and
\begin{equation}
 V'(u)=H(u,V(u)),\quad 
 V(u_0)=v(0),\quad
 t(u)=\int_{u_0}^u\frac{d\xi}{f(\xi,V(\xi))}
  \label{eq:DDDEg-pro12}
\end{equation}
holds.
\end{lem}
{\bf{Proof.}}
By the derivative of the inverse function, $t'(u)=1/\dot u=1/f(u,V(u))$. 
Moreover,
\begin{equation}
V'(u)=\frac{\dot v(t(u))}{\dot u(t(u))}=\frac{\alpha g(u)-\beta V(u)}{f(u,V(u))}=:H(u,V(u)), \quad 
V(u_0)=v(0)
  \label{eq:DDDEg-pro13}
\end{equation}
holds.
This completes the proof of Lemma \ref{lem:DDDEg-pro2}.
\qed
\\

From \eqref{eq:DDDEg-pro9} and \eqref{eq:DDDEg-pro12}, since $t(u)\to T_{\max}$ as $u\to U^*-0$, it follows immediately that
\begin{equation}
T_{\max}=\int_{u_0}^{U^*}\frac{d\xi}{f(\xi,V(\xi))}
\label{eq:DDDEg-pro14}
\end{equation}
holds.
With this preparation, we can prove the following lemma on the subsolution condition.

\begin{lem}[Subsolution condition]
\label{lem:DDDEg-pro3}
Let $u_1\in[u_0,U^*)$, and suppose that there exists a function $\Phi\in C^1([u_1,U^*);(0,v_{\max}))$ such that the following satisfies
\begin{equation}
 \Phi'(u)\le H(u,\Phi(u)) \quad(u_1\le u<U^*),\quad \Phi(u_1)\le v_{\min}.
 \label{eq:DDDEg-pro15}
\end{equation}
Then
\begin{equation}
 V(u)\ge\Phi(u),\quad 
 f(u,V(u))\ge f(u,\Phi(u)) \quad (u_1\le u<U^*)
 \label{eq:DDDEg-pro16}
\end{equation}
holds.
\end{lem}
{\bf{Proof.}}
By Lemma \ref{lem:DDDEg-pro1},
\begin{equation}
V(u_1) = v(t(u_1)) \ge v_{\min} \ge \Phi(u_1).
\label{eq:DDDEg-pro17}
\end{equation}
Since $H\in C^1(\Omega_+)$ satisfies the assumptions of Proposition \ref{prop:DDDEg-pre3}, we obtain $V\ge\Phi$ (see Remark \ref{rem:DDDEg-pre1}). 
The second inequality follows from $f_v\ge0$.
\qed
\\

Since $dt=du/f(u,V(u))$, setting $t_1=t(u_1)$ we have
\begin{equation}
T_{\max} - t_1 
=\int_{u_1}^{U^*}\frac{du}{f(u,V(u))}
\le \int_{u_1}^{U^*}\frac{du}{f(u, \Phi(u))}.
\label{eq:DDDEg-pro18}
\end{equation}
The case $U^*=\infty$ corresponds to blow-up in the true sense, while the case $U^*<\infty$ with $v_{\max}<\infty$ involves a singularity of the type $f\to\infty$ associated with $v\to v_{\max}$.

In constructing a subsolution $\Phi$, the practical point is how to guarantee $\Phi(u_1)\le V(u_1)$. 
Since $V(u)\ge v_{\min}$ for every $u\ge u_0$, if $\Phi(u_1)\le v_{\min}$ then
\begin{equation}
\Phi(u_1)\le v_{\min}\le V(u_1)
\label{eq:DDDEg-pro19}
\end{equation}
holds.
Hence the key point of the construction is to take $u_1$ sufficiently large so that $\Phi(u_1)\le v_{\min}$ and the differential inequality hold.

\subsection{Estimate of the blow-up time}
\label{sub:DDDEg-pro3}
In this subsection we show $T_{\max}<+\infty$. 
We can give the following upper bound for $T_{\max}$.

\begin{lem}
\label{lem:DDDEg-pro4}
Assume (I), (II), and (III-2). 
Set
\begin{equation}
c_1 := \min_{u\in[u_0,u_1]}f(u,v_{\min}).
\label{eq:DDDEg-pro20}
\end{equation}
Then $c_1>0$ holds, and the following holds
\begin{equation}
T_{\max} \le \dfrac{u_1-u_0}{c_1} + \int_{u_1}^\infty \dfrac{du}{f(u,\Phi(u))}
<\infty.
\label{eq:DDDEg-pro21}
\end{equation}
\end{lem}
{\bf{Proof.}}
Since $f$ is continuous and positive and $[u_0,u_1]$ is compact, $c_1>0$. 
By $f_v\ge0$ and Lemma \ref{lem:DDDEg-pro1}, as long as $u(t)\le u_1$,
\[
\dot u=f(u,v)\ge f(u,v_{\min})\ge c_1,
\]
and hence $u(t)\ge u_0+c_1t$. 
We now distinguish the following two cases (A) and (B).
\begin{enumerate}
\item[(A)] 
Suppose $u(t)<u_1$ for all $t\in[0,T_{\max})$. 
Then $u_0+c_1t\le u(t)<u_1$ gives $T_{\max}\le(u_1-u_0)/c_1$, and \eqref{eq:DDDEg-pro21} follows.
\item[(B)]
Suppose $u(t_1)=u_1$ for some $t_1\in[0,T_{\max})$. 
Since $u\le u_1$ on $[0,t_1]$, the same argument gives $u_1\ge u_0+c_1t_1$, that is, $t_1\le(u_1-u_0)/c_1$. 
In this case $U^*>u_1$, and by Lemma \ref{lem:DDDEg-pro2},
\begin{equation}
T_{\max}=t_1+\int_{u_1}^{U^*}\frac{d\xi}{f(\xi,V(\xi))}
\label{eq:DDDEg-pro22}
\end{equation}
holds.
Since $\Phi$ in assumption (III-2) is defined on $[u_1,\infty)$, we may restrict it to $[u_1,U^*)$. 
By Lemma \ref{lem:DDDEg-pro3}, $f(\xi,V(\xi))\ge f(\xi,\Phi(\xi))$, and therefore
\begin{equation}
 T_{\max}\le t_1+\int_{u_1}^{U^*}\frac{d\xi}{f(\xi,\Phi(\xi))}\le\frac{u_1-u_0}{c_1}+\int_{u_1}^{\infty}\frac{d\xi}{f(\xi,\Phi(\xi))}
\label{eq:DDDEg-pro23}
\end{equation}
holds.
By (III-2), the right-hand side is finite.
\end{enumerate}
This completes the proof of Lemma \ref{lem:DDDEg-pro4}.
\qed
\\

\subsection{Proof that $U^*=\infty$}
\label{sub:DDDEg-pro4}
Next we show $U^*=\infty$. 
By Lemma \ref{lem:DDDEg-pro4}, $T_{\max}<\infty$ holds. 
When $v_{\max}<\infty$, $T_{\max}<\infty$ alone leaves open the possibility that $f\to\infty$ as $v\to v_{\max}$ while $u$ remains bounded. 
The purpose of this subsection is to show that assumption (V) excludes this possibility.

\begin{lem}
\label{lem:DDDEg-pro5}
Let $v_{\max}=\infty$. 
If $T_{\max}<\infty$, then $U^*=\infty$.
\end{lem}
{\bf{Proof.}}
We argue by contradiction. 
Suppose $U^*<\infty$, and set
\begin{equation}
M := \max_{[u_0,U^*]}g <\infty.
\label{eq:DDDEg-pro24}
\end{equation}
From $g(u(t))\le M$ and the right-hand side of \eqref{eq:DDDEg-pro4}, we obtain
\begin{equation}
 v(t)\le e^{-\beta t}v(0)+\frac{\alpha M}{\beta}(1-e^{-\beta t}) 
 \le \max\left\{v(0),\frac{\alpha M}{\beta} \right\}
 =: \tilde{V} <\infty.
 \label{eq:DDDEg-pro25}
\end{equation}
Hence, we have $(u(t),v(t))\in K:=[u_0,U^*]\times[v_{\min},\widetilde V]$. 
By $v_{\max}=\infty$, $K$ is a compact subset of $\Omega_+$. 
However, we assumed $T_{\max}<\infty$, and if $T_{\max}<\infty$, then $(u(t),v(t))$ cannot remain in any compact subset of $\Omega_+$ by the continuation theorem. 
This is a contradiction. 
This completes the proof of Lemma~\ref{lem:DDDEg-pro5}.
\qed
\\

\begin{lem}
\label{lem:DDDEg-pro6}
Let $v_{\max}<\infty$, and suppose $T_{\max}<\infty$ and $U^*<\infty$. 
Then $v$ extends to a Lipschitz continuous function on $[0,T_{\max}]$, and we obtain
\begin{align}
& v(t)\to v_{\max} \quad {\rm{as}} \quad t\to T_{\max}-0,
\label{eq:DDDEg-pro26}
\\
& v_{\max}-v(t)\le L_0(T_{\max}-t) \quad (0\le t<T_{\max}),
\label{eq:DDDEg-pro27}
\end{align}
where $M=\max_{[u_0,U^*]}g$ and $L_0:=\alpha M+\beta v_{\max}$.
\end{lem}
{\bf{Proof.}}
By the existence of the constant $M$ and \eqref{eq:DDDEg-int5},
\begin{equation}
|\dot v| = |\alpha g(u)-\beta v| 
\le \alpha M+\beta v_{\max}
=L_0
\label{eq:DDDEg-pro28}
\end{equation}
holds.
Thus $v$ is Lipschitz continuous on $[0,T_{\max})$ with Lipschitz constant $L_0$. 
Hence $v$ is uniformly continuous and extends continuously to $t=T_{\max}$, and the limit
\begin{equation}
L_v := \lim_{t\to T_{\max}^{-}} v(t) \in [v_{\min},v_{\max}]
\label{eq:DDDEg-pro29}
\end{equation}
exists. 
It remains to show $L_v=v_{\max}$. 
Suppose $L_v<v_{\max}$. 
Then $(u(t),v(t))\to(U^*,L_v)$, which is an interior point of $\Omega_+$. 
Since $f$ and $g$ are of class $C^1$ on $\Omega_+$, the local existence theorem applies with this point as a new initial value, and the solution can be extended beyond $T_{\max}$. 
This contradicts the maximality of $T_{\max}$ and the continuation theorem. 
Therefore $L_v=v_{\max}$, and 
\begin{equation}
v_{\max}-v(t)=\lim_{t\to T_{\max}^{-}} v(s) - v(t) \le L_0(T_{\max}-t)
\label{eq:DDDEg-pro30}
\end{equation}
holds.
\qed
\\

\begin{lem}
\label{lem:DDDEg-pro7}
Let $v_{\max}<\infty$ and assume (V). 
If $T_{\max}<\infty$, then $U^*=\infty$.
\end{lem}
{\bf{Proof.}}
We argue by contradiction. 
Suppose $U^*<\infty$. 
By Lemma \ref{lem:DDDEg-pro6}, there exists $t_2\in[0,T_{\max})$ such that $v(t)>v_{\max}-\delta_0$ for $t\ge t_2$. 
By (V) and $u(t)\ge u_0$, for $t_2\le t<T_{\max}$, we obtain 
\begin{equation}
\dot u=f(u,v)\ge\kappa\,(v_{\max}-v(t))^{-q}\ge\kappa L_0^{-q}(T_{\max}-t)^{-q}.
\label{eq:DDDEg-pro31}
\end{equation}
Here the last inequality follows from $0<v_{\max}-v(t)\le L_0(T_{\max}-t)$ and the monotone decrease of $s\mapsto s^{-q}$. 
Integrating from $t_2$ to $t$, we obtain
\begin{equation}
u(t) \ge u(t_2)+\kappa L_0^{-q}\int_{t_2}^t(T_{\max}-s)^{-q}ds
\to \infty \quad
{\rm{as}} \quad t\to T_{\max}^{-},
\label{eq:DDDEg-pro32}
\end{equation}
which contradicts the assumption $U^*<\infty$. 
This completes the proof of Lemma \ref{lem:DDDEg-pro7}.
\qed
\\

\subsection{Proof of Theorem \ref{th:DDDEg-mr1} and Corollary \ref{cor:DDDEg-mr1}}
\label{sub:DDDEg-pro5}
By Lemma \ref{lem:DDDEg-pro4}, $T_{\max}<\infty$. 
Applying Lemma \ref{lem:DDDEg-pro5} when $v_{\max}=\infty$, and Lemma \ref{lem:DDDEg-pro7} under assumption (V) when $v_{\max}<\infty$, we obtain $U^*=\infty$, that is,
\[
  \lim_{t\to T_{\max}-0}u(t)=\infty.
\]
Hence \eqref{eq:DDDEg-pre2} holds with $T=T_{\max}\in(0,\infty)$. 
Furthermore, if (IV) holds, then the conditions of Definition \ref{def:DDDEg-mr1} are satisfied, so this blow-up is memory-induced blow-up. 
This completes the proof of Theorem \ref{th:DDDEg-mr1}. 
\qed
\\

For Corollary \ref{cor:DDDEg-mr1}, in assumption (III-2) take $\Phi(u)\equiv v_{\min}$ and $u_1:=u_0$. 
Then the range of $\Phi$ is $v_{\min}\in(0,v_{\max})$, and $\Phi(u_1)=v_{\min}$. 
By the monotonicity of $g$, for $u\ge u_0$,
\begin{equation}
\alpha g(u)-\beta v_{\min} \ge \alpha g(u_0)-\beta v_{\min} \ge 0
\label{eq:DDDEg-pro33}
\end{equation}
so that $\Phi'=0\le H(u,\Phi(u))$. 
Moreover, \eqref{eq:DDDEg-mr3} becomes $\int_{u_0}^{\infty}du/f(u,v_{\min})<\infty$, which is exactly (III-1). 
Hence the assumptions of Theorem \ref{th:DDDEg-mr1} are satisfied, and the proof of Corollary \ref{cor:DDDEg-mr1} is complete.
\qed
\\

\subsection{Proof of Theorem \ref{th:DDDEg-mr2}}
\label{sub:DDDEg-pro6}
Since $v_{\max}=\infty$, if $T_{\max}<\infty$ then $U^*=\infty$ by Lemma \ref{lem:DDDEg-pro5}. 
Define a function $F=F(u)$ by
\begin{equation}
F(u) := f \left(u,\max\left\{v(0),\frac{\alpha}{\beta}g(u)\right\} \right).
\label{eq:DDDEg-pro34}
\end{equation}
Substituting $t=t(u)$ into \eqref{eq:DDDEg-pro6} shown in Lemma \ref{lem:DDDEg-pro1} and using $g(u(t(u)))=g(u)$, we obtain
\[
 V(u) = v(t(u)) \le \max \left\{v(0), \dfrac{\alpha}{\beta}g(u) \right\} \quad 
 (u_0\le u<U^*).
\]
By $f_v\ge0$, $f(u,V(u))\le F(u)$. 
Therefore, by Lemma \ref{lem:DDDEg-pro2},
\begin{equation}
 T_{\max}=\int_{u_0}^{\infty}\frac{d\xi}{f(\xi,V(\xi))}\ge\int_{u_0}^{\infty}\frac{d\xi}{F(\xi)}
\label{eq:DDDEg-pro36}
\end{equation}
holds.
\qed
\\

\subsection{Proof of Theorem \ref{th:DDDEg-mr3}}
\label{sub:DDDEg-pro7}
We assume $T_{\max}<\infty$ and derive a contradiction. 
Using $H$ in \eqref{eq:DDDEg-pro11}, the second condition in \eqref{eq:DDDEg-mr9} can be written as $\Psi'(u)\ge H(u,\Psi(u))$ for $u\ge u_0$.

First, we show $V\le\Psi$ by the comparison theorem (Proposition \ref{prop:DDDEg-pre3}). 
On $[u_0,U^*)$, both $\Psi$ and $V$ are of class $C^1$ with values in $(0,v_{\max})$, and
\begin{equation}
 \Psi'(u)\ge H(u,x(u)), \quad 
 V'(u)=H(u,y(u)) \quad(u_0\le u<U^*)
 \label{eq:DDDEg-pro37}
\end{equation}
together with $\Psi(u_0)\ge v(0)=V(u_0)$. 
Here the equation for $V$ is due to Lemma \ref{lem:DDDEg-pro2}. 
Since $H\in C^1(\Omega_+)$ satisfies the assumptions of Proposition \ref{prop:DDDEg-pre3}, applying it with $(\Phi,\Psi)=(V,\Psi)$ (see Remark \ref{rem:DDDEg-pre1}) gives
\begin{equation}
V(u)\le\Psi(u)\qquad(u_0\le u<U^*).
\label{eq:DDDEg-pro38}
\end{equation}
Moreover, by $f_v\ge0$, $f(u,V(u))\le f(u,\Psi(u))$ holds.

Next, we show $U^*=\infty$. 
Suppose $U^*<\infty$. 
Since $\Psi$ is continuous on $[u_0,U^*]$ and takes values in $(0,v_{\max})$,
\begin{equation}
 \overline{\Psi}:=\max_{u\in[u_0,U^*]}\Psi(u)<v_{\max}
 \label{eq:DDDEg-pro39}
\end{equation}
holds.
Then, by Lemma~\ref{lem:DDDEg-pro1}, \eqref{eq:DDDEg-pro38}, and \eqref{eq:DDDEg-pro39}, for $t\in[0,T_{\max})$,
\begin{equation}
 v_{\min}\le v(t)=V(u(t))\le\Psi(u(t))\le\overline{\Psi}
 \label{eq:DDDEg-pro40}
\end{equation}
holds.
Hence $(u(t),v(t))$ remains in the compact set $[u_0,U^*]\times[v_{\min},\overline{\Psi}]\subset\Omega_+$. 
Since $T_{\max}<\infty$, this contradicts the continuation theorem. 
Therefore,  we obtain $U^*=\infty$.

Consequently, combining \eqref{eq:DDDEg-pro14} of Lemma \ref{lem:DDDEg-pro2} with \eqref{eq:DDDEg-pro38}, we obtain
\[
  T_{\max}=\int_{u_0}^{\infty}\frac{d\xi}{f(\xi,V(\xi))}\ge\int_{u_0}^{\infty}\frac{d\xi}{f(\xi,\Psi(\xi))}=\infty,
\]
which contradicts $T_{\max}<\infty$. 
Hence we have $T_{\max}=\infty$.
\qed\\

\section{Examples}
\label{sec:DDDEg-ex}
In this section, we present several examples in which Theorem \ref{th:DDDEg-mr1} and Corollary \ref{cor:DDDEg-mr1} are used to decide whether memory-induced blow-up occurs. 
The key point is how to find the subsolution $\Phi$ in (III-2). 
We first give a ``recipe'' (Proposition \ref{prop:DDDEg-ex1}) that constructs $\Phi$ mechanically when $f$ is of separable form, and then use it to examine several examples, including (F1) and (F2).

\subsection{Construction of the subsolution $\Phi$}
\label{sub:DDDEg-ex1}
Multiplying the inequality in (III-2), \eqref{eq:DDDEg-mr2}, by $f>0$, it is equivalent to
\begin{equation}
f(u,\Phi(u))\Phi'(u)+\beta\Phi(u) \le \alpha g(u)
\quad (u\ge u_1).
\label{eq:DDDEg-ex1}
\end{equation}
The term $\alpha g$ on the right-hand side is the amount of memory generated per unit time by the present state.
On the left-hand side, the second term $\beta\Phi$ is the loss due to forgetting, and the first term is the amount needed for $\Phi$ itself to increase. 
We therefore split $\alpha g$ into $\kappa g+(\alpha-\kappa)g$ ($0<\kappa<\alpha$) and reduce the problem to finding $\Phi$ satisfying
\begin{equation}
f(u,\Phi)\Phi' = \kappa g(u), \quad 
\beta\Phi(u) \le (\alpha-\kappa) g(u).
\label{eq:DDDEg-ex2}
\end{equation}
The first relation determines $\Phi$ by solving for it, and the second is used to check that $\Phi$ grows sufficiently more slowly than $g$.

We make the following setting. 
Let $v_{\max}\in(0,\infty]$. 
In addition to (I) and (II), assume that $f$ is separable:
\begin{equation}
f(u,v) = a(u)b(v),
\label{eq:DDDEg-ex3}
\end{equation}
where $a=a(u)$ and $b=b(v)$ satisfy
\begin{equation}
a>0, \quad a\in C^1((0,\infty)), \quad 
b>0, \quad b'\ge 0, \quad b\in C^1((0,v_{\max})).
\label{eq:DDDEg-ex4}
\end{equation}
Then the first relation in \eqref{eq:DDDEg-ex2} becomes
\begin{equation}
a(u)b(\Phi)\Phi' = \kappa g(u)
\label{eq:DDDEg-ex5}
\end{equation}
and, thanks to the assumption on $a$ in \eqref{eq:DDDEg-ex4}, the variables $u$ and $\Phi$ separate completely:
\begin{equation}
b(\Phi)\Phi' = \kappa \dfrac{g(u)}{a(u)}.
\label{eq:DDDEg-ex6}
\end{equation}
Furthermore, for $v_{\min}$ defined in \eqref{eq:DDDEg-mr1}, set
\begin{equation}
B(v) := \int_{v_{\min}}^{v} b(s)\,ds \quad (v_{\min}\le v<v_{\max}),\quad
\int_{v_{\min}}^{v_{\max}} b(s)\,ds = \infty.
\label{eq:DDDEg-ex7}
\end{equation}
Integrating \eqref{eq:DDDEg-ex5} from $u_1$ to $u$ gives
\[
\int_{u_1}^{u} b(\Phi(s)) \dfrac{d\Phi}{ds} \,ds = \kappa \int_{u_1}^{u} \dfrac{g(s)}{a(s)} ds.
\]
Substituting $\xi=\Phi(s)$ on the left-hand side, with $\xi=\Phi(u_1)=v_{\min}$ at $s=u_1$ and $\xi=\Phi(u)$ at $s=u$, the definition \eqref{eq:DDDEg-ex6} of $B$ yields
\[
\int_{\Phi(u_1)}^{\Phi(u)} b(\xi) d\xi 
= \int_{v_{\min}}^{\Phi(u)} b(\xi) d\xi
= B(\Phi(u))
\]
and hence
\begin{equation}
B(\Phi(u)) = \kappa \int_{u_1}^{u} \frac{g(s)}{a(s)} \,ds
\label{eq:DDDEg-ex8}
\end{equation}
holds.
The second condition in \eqref{eq:DDDEg-ex6} holds automatically when $v_{\max}=\infty$, since $b$ is positive and nondecreasing. 
When $v_{\max}<\infty$, it holds, for example, for $b(v)=(v_{\max}-v)^{-q}$ ($q\ge1$), which corresponds to $q\ge1$ in assumption (V). 
In this case $B:[v_{\min},v_{\max})\to[0,\infty)$ is a strictly increasing $C^1$ bijection. 
Therefore, for $\kappa\in(0,\alpha)$ and $u_1\ge u_0$, we define
\begin{equation}
\Phi_{u_1}(u) := B^{-1} \left(\kappa\int_{u_1}^{u} \dfrac{g(s)}{a(s)} \,ds \right)
\quad (u\ge u_1).
\label{eq:DDDEg-ex9}
\end{equation}
Under this setting, the following proposition holds.

\begin{prop}
\label{prop:DDDEg-ex1}
Under the above setting, $\Phi_{u_1}\in C^1([u_1,\infty);[v_{\min},v_{\max}))$, and the following hold.
\begin{enumerate}
\item[(i)] $\Phi_{u_1}(u_1)=v_{\min}$ and $f(u,\Phi_{u_1}(u))\Phi_{u_1}'(u)=\kappa g(u)$ for $u\ge u_1$.
\item[(ii)] $\Phi_{u_1}$ satisfies \eqref{eq:DDDEg-mr2} if and only if
\begin{equation}
\beta \Phi_{u_1}(u) \le(\alpha-\kappa) g(u) \quad(u\ge u_1).
 \label{eq:DDDEg-ex10}
\end{equation}
\item[(iii)] $\displaystyle\int_{u_1}^\infty\frac{du}{f(u,\Phi_{u_1}(u))}=\frac{1}{\kappa}\int_{u_1}^\infty\frac{\Phi_{u_1}'(u)}{g(u)}\,du$.
\item[(iv)] 
If $u_0\le u_1\le u_1'$, then $\Phi_{u_1'}\le\Phi_{u_1}$ on $[u_1',\infty)$. 
In particular, if
\begin{equation}
\limsup_{u\to\infty} \dfrac{\beta\Phi_{u_0}(u)}{g(u)} < \alpha - \kappa
\label{eq:DDDEg-ex11}
\end{equation}
then $\Phi_{u_1}$ satisfies \eqref{eq:DDDEg-ex9} for $u_1$ sufficiently large.
\end{enumerate}
Consequently, if \eqref{eq:DDDEg-ex9} holds and the integral in (iii) is finite, then $\Phi=\Phi_{u_1}$ satisfies all the conditions of (III-2).
\end{prop}
{\bf{Proof.}} 
Since $g\ge g(u_0)>0$ for $u\ge u_0$, the function $g/a$ is continuous and positive, so the expression $\kappa\int_{u_1}^{u}(g/a)\,ds$ in the parentheses of \eqref{eq:DDDEg-ex8} is a nonnegative, increasing $C^1$ function of $u$. 
Since $B^{-1}$ is of class $C^1$ on $[0,\infty)$, $\Phi_{u_1}\in C^1$.
Moreover $B^{-1}(0)=v_{\min}$ gives $\Phi_{u_1}(u_1)=v_{\min}$, and the range is contained in $[v_{\min},v_{\max})\subset(0,v_{\max})$.
\begin{enumerate}
\item[(i)] 
Differentiating both sides of \eqref{eq:DDDEg-ex7} in $u$ and using \eqref{eq:DDDEg-ex6}, we get
\begin{equation}
b(\Phi_{u_1}(u))\Phi_{u_1}'(u) = \kappa \frac{g(u)}{a(u)}.
\label{eq:DDDEg-ex12}
\end{equation}
Multiplying both sides by $a(u)$ gives $a(u)b(\Phi_{u_1})\Phi_{u_1}'=\kappa g$, and \eqref{eq:DDDEg-ex3} yields (i).
\item[(ii)]
By (i), the inequality \eqref{eq:DDDEg-ex1}, which is equivalent to \eqref{eq:DDDEg-mr2}, becomes
\begin{equation}
\kappa g(u) + \beta \Phi_{u_1}(u) \le \alpha g(u),
\label{eq:DDDEg-ex13}
\end{equation}
which gives (ii).
\item[(iii)]
By (i), $f(u,\Phi_{u_1}(u))\Phi_{u_1}'(u)=\kappa g(u)$, so
\begin{equation}
\dfrac{1}{f(u, \Phi_{u_1}(u))} = \dfrac{\Phi_{u_1}'(u)}{\kappa g(u)}.
\label{eq:DDDEg-ex14}
\end{equation}
Integrating from $u_1$ to $\infty$ gives (iii). 
This means that, instead of integrating the original $f$ directly, it suffices to integrate the ratio of the derivative of the subsolution to $g$, which makes the convergence of the integral easy to discuss.
\item[(iv)]
For $u\ge u_1'$,
\[
\int_{u_1'}^{u} \dfrac{g}{a} \,ds \le \int_{u_1}^{u} \dfrac{g}{a}\,ds
\]
and since $B^{-1}$ is increasing, $\Phi_{u_1'}(u)\le\Phi_{u_1}(u)$. 
By \eqref{eq:DDDEg-ex10}, there exists $u_*\ge u_0$ such that $\beta\Phi_{u_0}(u)\le(\alpha-\kappa)g(u)$ for $u\ge u_*$. 
If $u_1\ge u_*$, then $\beta\Phi_{u_1}(u)\le\beta\Phi_{u_0}(u)\le(\alpha-\kappa)g(u)$ for $u\ge u_1$. 
This is a powerful property: even if it is difficult to find a subsolution satisfying the inequality for all $u\ge u_0$ from the start, a subsolution built by this recipe always satisfies the condition if we start from a sufficiently large $u_1$.
\end{enumerate}
This completes the proof of Proposition \ref{prop:DDDEg-ex1}.
\qed
\\
 
Based on the above, we summarize the verification procedure step by step. 
From the next subsection on, we present several applications of the criteria for memory-induced blow-up established in Section \ref{sec:DDDEg-mr}, following these steps.

\begin{quote}
\begin{enumerate}
\item[\textbf{Step 0}]: 
Check $v_{\min}$ and $v_{\max}$, and check that $f$ and $g$ satisfy assumptions (I) and (II). 
Then check assumption (III-1) to try Corollary \ref{cor:DDDEg-mr1}.
\item[\textbf{Step 1}]: 
Write $f$ in the form $a(u)b(v)$ and find $a$, $b$, and $B$.
\item[\textbf{Step 2}]: 
Examine the asymptotic behavior of $\Phi_{u_0}$ and check \eqref{eq:DDDEg-ex11}. 
If $\Phi_{u_0}$ grows more slowly than $g$, the left-hand side is $0$, and any $\kappa\in(0,\alpha)$ may be used.
\item[\textbf{Step 3}]: 
Compute or estimate the integral in condition (iii) of Proposition~\ref{prop:DDDEg-ex1}. 
The condition under which this integral is finite is the blow-up condition.
\item[\textbf{Step 4}]: 
If $v_{\max}<\infty$, check (V). 
If (IV) holds, the blow-up is memory-induced.
\item[\textbf{Step 5}]: 
Compare with the necessary condition \eqref{eq:DDDEg-mr7} of Theorem \ref{th:DDDEg-mr2} and see whether the conditions coincide.
\end{enumerate}
\end{quote}

\begin{rem}
\label{rem:DDDEg-ex1}
Even if $f$ is not separable, if there exist $\bar f:=a(u)b(v)$ and a constant $c>0$ such that $c\bar f\le f\le\bar f$, then $\Phi_{u_1}$ constructed from $\bar f$ can be used as it is. 
Indeed, since $\Phi_{u_1}'>0$, we have $f\Phi_{u_1}'\le\bar f\Phi_{u_1}'=\kappa g$, so the estimate of the left-hand side of \eqref{eq:DDDEg-ex1} is unchanged. 
Moreover, $1/f(u,\Phi_{u_1})\le c^{-1}/\bar f(u,\Phi_{u_1})$, so finiteness of the integral is also preserved.
\end{rem}

\begin{rem}
\label{rem:DDDEg-ex2}
If the integral of $b$ from $v_{\min}$ to $v_{\max}$ in \eqref{eq:DDDEg-ex6} is finite, the construction of $\Phi$ in this paper does not apply, and a different device is needed.
\end{rem}

\subsection{Some examples}
\label{sub:DDDEg-ex2}
Let $\alpha$, $\beta>0$, and let $\phi\in C_b$ satisfy $g(u_0)>0$ (and $v(0)<v_{\max}$ if $v_{\max}<\infty$). 
Fix $\kappa\in(0,\alpha)$.

\begin{ex}[(F1): $f(u,v)=v^{p}$ ($p>1$), $g(u)=u$]
\label{ex:DDDEg-ex1}
Here $v_{\max}=\infty$, and $f$ is of class $C^1$ on $\Omega_+$ with
\[
f(u,v)=v^{p}>0, \quad 
f_{v}(u,v)=pv^{p-1} \ge 0 \quad (v>0).
\]
Moreover, $f$ extends continuously to $v=0$, and $f(u,0)=0$ is locally Lipschitz continuous, so (I) holds. 
Clearly $g(u)=u$ satisfies (II). 
However, $f(u,v_{\min})=v_{\min}^{p}$ and
\[
\int_{u_0}^{\infty} \frac{du}{v_{\min}^p} = \left[ \frac{u}{v_{\min}^p} \right]_{u_0}^{\infty} = \infty
\]
so (III-1) does not hold. 
Hence the test by Corollary \ref{cor:DDDEg-mr1} is not available, and we need to construct a dynamic subsolution $\Phi_{u_1}(u)$.

Since $f$ is separable, we may take $a\equiv1$ and $b(v)=v^{p}$. 
By \eqref{eq:DDDEg-ex7}, we obtain 
\[
B(v) 
= \int_{v_{\min}}^{v} s^p \,ds 
= \dfrac{v^{p+1} - v_{\min}^{p+1}}{p+1}
\]
and 
\begin{equation}
\int_{u_1}^{u} \frac{g(s)}{a(s)} \,ds 
= \int_{u_1}^{u} \frac{s}{1} \,ds 
= \frac{u^2 - u_1^2}{2}.
\label{eq:DDDEg-ex16}
\end{equation}
Hence \eqref{eq:DDDEg-ex8} gives
\[
\frac{\Phi_{u_1}(u)^{p+1} - v_{\min}^{p+1}}{p+1} = \kappa \frac{u^2 - u_1^2}{2}.
\]
Solving for $\Phi_{u_1}(u)$ (that is, applying $B^{-1}$), by \eqref{eq:DDDEg-ex9} we obtain
\begin{equation}
\Phi_{u_1}(u) = \left(v_{\min}^{p+1}+\dfrac{(p+1)\kappa}{2}(u^2-u_1^2) \right)^{1/(p+1)}.
 \label{eq:DDDEg-ex15}
\end{equation}
Since $p>1$, we have $1/(p+1)<1$, so using $(A+B)^{\gamma}\le A^{\gamma}+B^{\gamma}$ for $\gamma\in(0,1)$ and $A,B\ge0$, we get
\[
\Phi_{u_1}(u) \le v_{\min} + C u^{\gamma}, \quad \gamma := \dfrac{2}{p+1}\in (0,1).
\]
To check \eqref{eq:DDDEg-ex9}: since $g(u)=u$ and $\gamma<1$, letting $u\to+\infty$,
\[
\dfrac{\beta \Phi_{u_1}(u)}{g(u)} 
\le \dfrac{\beta (v_{\min} + C u^\gamma)}{u} 
\to 0 
\quad {\rm{as}} \quad u \to \infty.
\]
Therefore, by the argument of Proposition \ref{prop:DDDEg-ex1} (iv), condition \eqref{eq:DDDEg-ex10} holds for sufficiently large $u_1$. 
This corresponds to Step 2 above.

Next, computing the integral in Proposition \ref{prop:DDDEg-ex1} (iii), with $C$ a positive constant,
\begin{align*}
\int_{u_1}^{\infty} \dfrac{du}{f(u, \Phi_{u_1}(u))} 
&= \dfrac{1}{\kappa} \int_{u_1}^{\infty} \dfrac{\Phi_{u_1}'(u)}{g(u)} \,du \\
&\sim \dfrac{1}{\kappa} \int_{u_1}^{\infty} \dfrac{Cu^{\gamma-1}}{u} \,du
\quad {\rm{as}} \quad u \to \infty \\
&= \dfrac{1}{\kappa} \int_{u_1}^{\infty} Cu^{\gamma-2} \,du \\
&= \dfrac{1}{\kappa} \int_{u_1}^{\infty} Cu^{-\frac{2p}{p+1}} \,du \\
&<+\infty
\end{align*}
which is finite since $\gamma-2<-1$. 
Since $v_{\max}=+\infty$, by Theorem \ref{th:DDDEg-mr1}, the solution $u=u(t)$ of the distributed delay differential equation \eqref{eq:DDDEg-int1}, \eqref{eq:DDDEg-int2} under (F1) blows up at some finite time $T_{\max}$.

Moreover, (IV) holds because $f(u,0)=0$, that is, $\dot U=0$ has a global solution. 
Hence, by Theorem \ref{th:DDDEg-mr1}, the solution blows up in finite time, and this blow-up is memory-induced. 
This conclusion agrees with \cite{DDDEb}.

We also compute the integral
\[
\int_{u_0}^{\infty} \frac{du}{f(u, \max\{v(0), \frac{\alpha}{\beta}u\})}.
\]
For sufficiently large $u$, the integrand is $(\beta\alpha^{-1})^{p}u^{-p}$, so this integral converges if and only if $p>1$.
\qed
\\
\end{ex}

More precisely, set $c:=\alpha/\beta$ and $u_*:=\max\{u_0,v(0)/c\}$. 
For $u\ge u_*$ we have $cu\ge v(0)$, so $\max\{v(0),cu\}=cu$, and
\[
 \int_{u_0}^{\infty} \dfrac{du}{f\left(u,\max\{v(0),cu\}\right)}
 =\int_{u_0}^{u_*}\frac{du}{\max\{v(0),cu\}^{p}}
 +c^{-p}\int_{u_*}^{\infty}u^{-p}\,du.
\]
The first term on the right-hand side is finite, and the second diverges to $+\infty$ when $p\le1$. 
Hence, by the contrapositive of Theorem \ref{th:DDDEg-mr2}, for $0<p\le1$ the solution exists globally for any $\alpha,\beta>0$ and any $\phi\in C_b$. 
Combined with the result of this example for $p>1$, we conclude that the solution of (F1) blows up in finite time if and only if $p>1$.

For $0<p\le1$, the solution of the non-delayed equation $\dot U=(cU)^{p}$ with $c$ fixed also exists globally. 
For $p>1$, the blow-up is caused by the fact that the memory $v$ grows following $u$ and $v^{p}$ acts superlinearly. 
For $p\le1$, the same accumulation of memory is not superlinear enough, and blow-up does not occur.
\qed
\\ 
 
 \begin{ex}[(F2): $f(u,v)=(1-v)^{-q}$ ($q\ge1$), $g(u)=u$]
 \label{ex:DDDEg-ex2}
Here $v_{\max}=1$, and $f$ and $g$ satisfy (I) and (II). 
However, $f(u,v_{\min})=(1-v_{\min})^{-q}$ and
 \[
 \int_{u_0}^{\infty} \frac{du}{(1 - v_{\min})^{-q}} = \infty
 \]
so (III-1) does not hold. 
The test by Corollary \ref{cor:DDDEg-mr1} is not available, and we need to construct a dynamic subsolution $\Phi_{u_1}(u)$.
 
Since $f$ is separable, we may take $a(u)\equiv1$ and $b(v)=(1-v)^{-q}$. 
By \eqref{eq:DDDEg-ex7}:
 \begin{itemize}
 \item For $q=1$,
 \[
 B(v) 
 = \int_{v_{min}}^{v} \frac{1}{1-s} \,ds 
 = \left[ -\log(1-s) \right]_{v_{\min}}^{v} 
 = \log\frac{1-v_{\min}}{1-v}
 \]
and \eqref{eq:DDDEg-ex16} and \eqref{eq:DDDEg-ex8} give
 \[
  \log\frac{1-v_{\min}}{1-\Phi_{u_1}(u)}=\kappa \dfrac{u^{2}-u_{1}^{2}}{2}.
 \]
Solving for $\Phi_{u_1}(u)$ (that is, applying $B^{-1}$), by \eqref{eq:DDDEg-ex9},
\[
\Phi_{u_1}(u) = 1 - (1-v_{\min})\exp\left( -\frac{\kappa}{2}(u^2 - u_1^2) \right)
\]
holds.
 \item For $q>1$,
 \[
 B(v) 
= \int_{v_{min}}^{v} (1-s)^{-q} ds 
= \left[ \frac{(1-s)^{-q+1}}{q-1} \right]_{v_{\min}}^{v} 
= \frac{1}{q-1} \left( (1-v)^{-q+1} - (1-v_{\min})^{-q+1} \right)
 \]
and \eqref{eq:DDDEg-ex16} and \eqref{eq:DDDEg-ex8} give
\[
 \frac{1}{q-1} \left( (1-\Phi_{u_1}(u))^{-q+1} - (1-v_{\min})^{-q+1} \right) =\kappa \dfrac{u^{2}-u_{1}^{2}}{2}
 \]
holds.
Solving for $\Phi_{u_1}(u)$ (that is, applying $B^{-1}$), by \eqref{eq:DDDEg-ex9}, we obtain the following:
\[
\Phi_{u_1}(u) = 1 - \left[ (1-v_{\min})^{1-q} + \frac{(q-1)\kappa}{2}(u^2 - u_1^2) \right]^{-\frac{1}{q-1}}. 
\]
 \end{itemize}
In both cases $\Phi_{u_1}<1$, so if $u_1\ge\max\{u_0,\beta/(\alpha-\kappa)\}$, then $\beta\Phi_{u_1}<\beta\le(\alpha-\kappa)u$ for $u\ge u_1$, and hence \eqref{eq:DDDEg-ex11} holds.

Next, computing the integral in Proposition \ref{prop:DDDEg-ex1} (iii):
\begin{itemize}
\item
For $q=1$, $f(u,\Phi_{u_1}(u))=(1-\Phi_{u_1}(u))^{-1}$, and from the form of the subsolution, $1-\Phi_{u_1}(u)=(1-v_{\min})e^{-\frac{\kappa}{2}(u^{2}-u_1^{2})}$, so
\[
\frac{1}{f(u, \Phi_{u_1}(u))} = (1-v_{\min})\exp\left( -\frac{\kappa}{2}(u^2 - u_1^2) \right).
\]
Its integral is of Gaussian type and converges to a finite value.
\item 
For $q>1$, the integrand satisfies $(1-\Phi_{u_1})^{q}\sim c\,u^{-2q/(q-1)}$, and since $2q/(q-1)>1$, the integral is finite.
\end{itemize}
Assumption (V) holds as it is with $\kappa'=1$ and the same $q$ (we write $\kappa'$ to avoid confusion with the $\kappa$ of the recipe), and (IV) holds since $f(u,0)=1$ ($\dot U=1$). 
Hence, by Theorem \ref{th:DDDEg-mr1}, the solution blows up in finite time, and the blow-up is memory-induced. In particular, the conclusion for $q=1$ agrees with \cite{DDDEb}.
\qed\\
 \end{ex}
 
Examples \ref{ex:DDDEg-ex1} and \ref{ex:DDDEg-ex2} show that our criteria reproduce the results of \cite{DDDEb}. 
In \cite{DDDEb}, the dynamics of the two-dimensional system \eqref{eq:DDDEg-int5}, including infinity, was clarified by the Poincar\'e-type compactification, but when $f$ contains an exponential function, the method of \cite{DDDEb} cannot be applied immediately. 
The following example shows that the criteria established in this paper also apply to cases whose dynamics cannot be clarified by \cite{DDDEb}.

\begin{ex}[$f(u,v)=ue^{v}$, $g(u)=u$]
\label{ex:DDDEg-ex3}
Here $v_{\max}=\infty$, and $f$ and $g$ satisfy assumptions (I) and (II) of Section \ref{sec:DDDEg-mr}. 
Since $f(u,v_{\min})=ue^{v_{\min}}$,
\[
\int_{u_0}^{\infty} \frac{du}{u e^{v_{\min}}} 
= \frac{1}{e^{v_{\min}}} \left[\dfrac{}{}\log u \dfrac{}{}\right]_{u_0}^{\infty} 
= \infty
\]
so the test by Corollary \ref{cor:DDDEg-mr1} is not available. 
As in the previous examples, we construct a dynamic subsolution $\Phi_{u_1}(u)$.

Taking $a(u)=u$ and $b(v)=e^{v}$, $f$ is separable, and
\[
B(v) = \int_{v_{\min}}^{v} e^s \,ds = e^v - e^{v_{\min}}, \quad
\int_{u_1}^{u} \frac{g(s)}{a(s)} ds = \int_{u_1}^{u} \frac{s}{s} ds = u - u_1
\]
holds.
By \eqref{eq:DDDEg-ex8},
\[
e^{\Phi_{u_1}(u)} - e^{v_{\min}} = \kappa (u - u_1)
\]
and solving for $\Phi_{u_1}(u)$ (that is, applying $B^{-1}$), by \eqref{eq:DDDEg-ex9}, we obtain
\[
\Phi_{u_1}(u) = \log\left( e^{v_{\min}} + \kappa(u - u_1) \right).
\]
Taking $u_1\ge1$ sufficiently large, we have
\[
 \Phi_{u_1}(u) \le \log(e^{v_{\min}}+\kappa) + \log u
 \quad (u\ge u_1\ge1)
\]
so that 
\[
\dfrac{\beta \Phi_{u_1}(u)}{g(u)} \le \dfrac{\beta \left( \log(e^{v_{\min}} + \kappa) + \log u \right)}{u} 
\to  0 
\quad {\rm{as}} \quad u \to \infty
\]
holds.
By Proposition \ref{prop:DDDEg-ex1} (iv), \eqref{eq:DDDEg-ex10} holds for sufficiently large $u_1$.

Computing the integral in Proposition \ref{prop:DDDEg-ex1} (iii), we have
\[
\int_{u_1}^{\infty} \frac{du}{f(u, \Phi_{u_1}(u))} = \int_{u_1}^{\infty} \frac{du}{u e^{\Phi_{u_1}(u)}} = \int_{u_1}^{\infty} \frac{du}{u \left( e^{v_{\min}} + \kappa(u - u_1) \right)}
<\infty.
\]
(IV) holds since $f(u,0)=u$ ($\dot U=U$) has a global solution. 
Therefore, by Theorem \ref{th:DDDEg-mr1}, the solution blows up in finite time, and the blow-up is memory-induced.
\qed\\
\end{ex}

\begin{ex}[Power type: $f(u,v)=u^{m}v^{p}$ ($m\in\mathbb{R}$, $p>0$), $g(u)=u$]\label{ex:DDDEg-ex4}
Here $v_{\max}=\infty$, and $f$ and $g$ satisfy (I) and (II). 
Moreover, (IV) holds since $f(u,0)=0$, that is, $\dot U=0$ has a global solution. 
Since $f(u,v_{\min})=u^{m}v_{\min}^{p}$,
\[
\dfrac{1}{f(u, v_{\min})} = \dfrac{1}{u^m v_{\min}^p} = v_{\min}^{-p} u^{-m}
\]
so for $m>1$ assumption (III-1) holds, and by Corollary \ref{cor:DDDEg-mr1}, memory-induced blow-up occurs.

We now examine this condition further by constructing a dynamic subsolution. 
Since $f$ is separable, take $a(u)=u^{m}$ and $b(v)=v^{p}$. 
Then we have 
\[
B(v) = \int_{v_{\min}}^{v} s^p \,ds 
= \dfrac{v^{p+1} - v_{\min}^{p+1}}{p+1}, \quad
\int_{u_{1}}^{u} \dfrac{g(s)}{a(s)}\, ds = \int_{u_1}^{u} s^{1-m} ds =:S(u)
\]
so \eqref{eq:DDDEg-ex8} gives
\[
\dfrac{\Phi_{u_{1}(u)}^{p+1} - v_{\min}^{p+1}}{p+1}=\kappa S(u)
\]
and solving for $\Phi_{u_1}(u)$ (that is, applying $B^{-1}$),
\[
\Phi_{u_1}(u) = \left( v_{\min}^{p+1} + (p+1)\kappa S(u) \right)^{\frac{1}{p+1}}.
\]
To estimate the growth of this subsolution, we look at the behavior of $S(u)$ and distinguish cases according to $m$:
\begin{itemize}
\item 
For $m<2$, $S(u)\le u^{2-m}/(2-m)$, so with constants $C$, $C'$, we obatain
\[
\Phi_{u_1}(u) 
\le \left( v_{\min}^{p+1} + C u^{2-m} \right)^{\frac{1}{p+1}} 
\le v_{\min} + C' u^{\gamma}
\] 
with $\gamma=(2-m)/(p+1)$.
Hence, if $\gamma<1$ (that is, $m+p>1$), then $\beta\Phi_{u_1}(u)/g(u)\to0$ as $u\to\infty$.
\item 
For $m=2$, $S(u)=\log(u/u_1)$, so $\Phi_{u_1}(u)=\mathcal{O}((\log u)^{1/(p+1)})$ and $\beta\Phi_{u_1}(u)/g(u)\to0$ as $u\to\infty$.
\item 
For $m>2$, $\Phi_{u_1}$ is bounded, and $\Phi_{u_1}/u\to0$ as $u\to\infty$.
\end{itemize}
Next, consider the integral in Proposition \ref{prop:DDDEg-ex1} (iii). 
For $m>1$, $\Phi_{u_1}(u)\ge v_{\min}>0$, so $1/f(u,\Phi_{u_1})\le v_{\min}^{-p}u^{-m}$, which is integrable. 
For $m\le1$ (so $m<2$), by the above, $\Phi_{u_1}(u)\ge cu^{\gamma}$ for sufficiently large $u$ and some $c>0$, hence $1/f(u, \Phi_{u_1}(u))=\mathcal{O}(u^{-(m+p\gamma)})$. 
The integral is finite if the exponent is larger than $1$, that is, if
\[
m+p\gamma=\dfrac{m+2p}{p+1}>1
\]
which is equivalent to $m+p>1$. 
Therefore, by Theorem \ref{th:DDDEg-mr1}, blow-up occurs and is memory-induced.
\qed
\\
\end{ex}

By using a dynamic subsolution, we obtain a blow-up condition in a range wider than $m>1$ in Corollary \ref{cor:DDDEg-mr1}. 
Moreover, in this case we can also prove that $m+p>1$ is necessary for finite-time blow-up. 
Indeed, for sufficiently large $u$,
\[
f\left(u,\max\{v(0),\frac{\alpha}{\beta}u\}\right)
=\frac{1}{u^m \left(\frac{\alpha}{\beta}u\right)^p} 
= \left(\frac{\beta}{\alpha}\right)^p u^{-(m+p)}
\]
so for $m+p>1$,
\[
\int_{u_0}^\infty \dfrac{du}{f\left(u,\max\{v(0),\frac{\alpha}{\beta}u\}\right)}<\infty
\]
whereas for $m+p\le1$ the integral diverges, and by Theorem \ref{th:DDDEg-mr2} the solution exists globally. 
Thus, in Example \ref{ex:DDDEg-ex4}, the sufficient and necessary conditions coincide, and we obtain that $m+p>1$ is a necessary and sufficient condition for finite-time blow-up.

\begin{ex}[$f(u,v)=(1-v)^{-q}$ ($q\ge1$), $g(u)=u/(1+u)$]
\label{ex:DDDEg-ex5}
Let $v_{\max}=1$ and $c:=\alpha/\beta$. 
By the same argument as in Example \ref{ex:DDDEg-ex2}, $f$ satisfies (I). 
On the other hand, $g(u)$ is of class $C^1$ on $[0,\infty)$ with
\[
g(u) \in [0, 1), \quad
g'(u)=\dfrac{1}{(1+u)^2} \ge 0, \quad
g(u_0) > 0
\]
so (II) holds. 
Let the initial function $\phi\in C_b$ satisfy
\begin{equation}
 c<1, \quad 
 v(0) \le cg(u_0).
 \label{eq:DDDEg-ex17}
\end{equation}
Define the function $\Psi=\Psi(u)$ by
\[
\Psi(u) := cg(u) = \dfrac{cu}{1+u}.
\]
Since $g$ is increasing with $g(u_0)>0$ and $g<1$,
\[
0 < \Psi(u_0) \le \Psi(u) < c < 1=v_{\max} \quad (u\ge u_0)
\]
and $\Psi\in C^1([u_0,\infty);(0,v_{\max}))$. 
By \eqref{eq:DDDEg-ex17}, $\Psi(u_0)=cg(u_0)\ge v(0)$. 
Moreover, $\alpha g-\beta\Psi=\alpha g-\beta cg=0$, so
\[
\frac{\alpha g(u)-\beta\Psi(u)}{f(u,\Psi(u))}=0\le cg'(u)=\Psi'(u).
\]
This is the case where condition \eqref{eq:DDDEg-mr10} of Remark \ref{rem:DDDEg-mr2} is satisfied.

Now $1/f(u,\Psi(u))=(1-cg(u))^{q}\ge(1-c)^{q}>0$, so
\[
 \int_{u_0}^{\infty}\frac{du}{f(u,\Psi(u))}
 \ge(1-c)^{q}\int_{u_0}^{\infty}du=\infty.
\]
Hence all the conditions \eqref{eq:DDDEg-mr9} of Theorem~\ref{th:DDDEg-mr3} hold, and under \eqref{eq:DDDEg-ex17} the solution exists globally ($T_{\max}=\infty$). 
Furthermore, by Lemma \ref{lem:DDDEg-pro1}, $v(t)\le\max\{v(0),cg(u(t))\}<1$, so $v$ never reaches $v_{\max}=1$.
\qed
\\
\end{ex}

Examples \ref{ex:DDDEg-ex2} and \ref{ex:DDDEg-ex5} show that, for the same $f=(1-v)^{-q}$, the conclusion is reversed solely by the boundedness of $g$. 
For $g(u)=u$, the memory $v$ reaches $v_{\max}=1$ and blow-up occurs. 
On the other hand, when $g$ is bounded and $c\sup g<1$, $v$ does not exceed $c\sup g<1$ and never reaches the singularity. 
Since $v_{\max}<\infty$, Theorem \ref{th:DDDEg-mr2} cannot be used, and Theorem \ref{th:DDDEg-mr3} is needed to prove global existence. 
The cases $c\ge1$ and $v(0)>cg(u_0)$ cannot be treated by the method of this paper; constructing criteria that also cover these cases is left for future work.

\begin{ex}[$f=u^{2}+v$, $g(u)=u$]
\label{ex:DDDEg-ex6}
Let $v_{\max}=\infty$. 
The function $f(u,v)=u^{2}+v$ is of class $C^1$ on $\Omega_+$ with
\[
 f(u,v)=u^{2}+v>0,\quad 
 f_v(u,v)=1
\]
extends continuously to $v=0$, and $f(u,0)=u^{2}$ is locally Lipschitz continuous. 
Hence (I) holds, and $g(u)=u$ satisfies (II).
For 
\[
v_{\min} = \min\left\{v(0),\dfrac{\alpha}{\beta}u_0 \right\}>0,
\]
we have $f(u,v_{\min})=u^{2}+v_{\min}$, so
\[
 \int_{u_0}^{\infty} \dfrac{du}{u^{2}+v_{\min}}
 = \dfrac{1}{\sqrt{v_{\min}}} \Bigl(\frac{\pi}{2}-\arctan\frac{u_0}{\sqrt{v_{\min}}}\Bigr)
 <\infty
\]
and assumption (III-1) holds. 
Hence Corollary \ref{cor:DDDEg-mr1} applies directly, and the solution $u=u(t)$ blows up in finite time.

On the other hand, the non-delayed equation obtained in the limit $c=\alpha/\beta\to0$ is
\[
 \dot U=f(U,0)=U^{2},\quad U(0)=u_0,
\]
whose solution is known to blow up in finite time, so assumption (IV) does not hold. 
Therefore, by Definition \ref{def:DDDEg-mr1}, this blow-up is not memory-induced.
\qed\\
\end{ex}

In Example~\ref{ex:DDDEg-ex6}, the memory term $v>0$ merely accelerates the blow-up through $\dot u=u^{2}+v>u^{2}$, and the cause of the blow-up is the ODE part $u^{2}$, which is unrelated to memory. 
In contrast, in Examples \ref{ex:DDDEg-ex1}--\ref{ex:DDDEg-ex4}, $f(u,0)$ does not cause blow-up, so the blow-up is due solely to the accumulation of memory. 
Thus ``finite-time blow-up'' and ``memory-induced blow-up'' are determined by separate tests, and the latter is distinguished by (IV).

\section{Concluding Remarks}
\label{sec:DDDEg-dis}
We discuss the validity of the definition of memory-induced blow-up (Definition \ref{def:DDDEg-mr1}) introduced in this paper. 
Compared with the non-delayed equation $\dot U=f(U,cg(U))$ with the parameter $c=\alpha/\beta$ fixed, the following property holds (Proposition \ref{prop:DDDEg-dis1}): if the delayed system blows up in finite time, then the corresponding non-delayed system with the same parameter must also blow up.

\begin{prop}[Comparison with the non-delayed equation with $c$ fixed]
\label{prop:DDDEg-dis1}
Let $v_{\max}=\infty$ and assume (I) and (II). 
Assume further that either $g$ is unbounded or $v(0)\le\alpha\beta^{-1}g(u_0)$. 
Set $c:=\alpha/\beta$. 
If the solution of the delayed system blows up in finite time, then
\begin{equation}
 \int_{u_0}^\infty\frac{du}{f(u,cg(u))}<\infty,
 \label{eq:DDDEg-dis1}
\end{equation}
that is, the solution of the non-delayed equation $\dot U=f(U,c\,g(U))$, $U(0)=u_0$, also blows up in finite time.
\end{prop}
{\bf{Proof.}}
Since $v_{\max}=\infty$ and (I), (II) are assumed, \eqref{eq:DDDEg-mr7} holds by Theorem \ref{th:DDDEg-mr2}. 
First consider the case $v(0)\le cg(u_0)$. 
Since $g$ is nondecreasing, $\max\{v(0),cg(u)\}=cg(u)$ for all $u\ge u_0$, and the claim follows.

Next consider the case where $g$ is unbounded. 
Since $g$ is nondecreasing and unbounded, there exists $u_*\ge u_0$ such that $cg(u_*)\ge v(0)$. 
For $u\ge u_*$ we have $\max\{v(0),cg(u)\}=cg(u)$, so
\begin{equation}
\int_{u_*}^\infty \dfrac{du}{f(u,cg(u))} < \infty
\label{eq:DDDEg-dis2}
\end{equation}
holds.
On the other hand, the function $h(u):=1/f(u,cg(u))$ is continuous on $[u_0,u_*]$, so
\begin{equation}
\int_{u_0}^{u_*} \dfrac{du}{f(u,cg(u))} < \infty
\label{eq:DDDEg-dis3}
\end{equation}
holds.
Hence the integral over the whole interval $[u_0,\infty)$ is finite. 
Finally, $F(u):=f(u,cg(u))$ is positive and locally Lipschitz continuous for $u\ge u_0$, so by Proposition \ref{prop:DDDEg-pre4}, finite-time blow-up of the solution of $\dot U=F(U)$ is equivalent to finiteness of the integral.
\qed
\\

Therefore, in order to capture the phenomenon in which ``the non-delayed system has a global solution, whereas only the delayed system blows up'' as an effect of the delay term (memory effect) alone, it is natural not to fix the parameter $c$ but to take as the reference the limit $c\to0$ (that is, $\dot U=f(U,0)$). 
In this sense, the formulation of Definition \ref{def:DDDEg-mr1} is a natural choice. 
We have not treated the case $v_{\max}<\infty$, nor the case where $g$ is bounded and $v(0)>\alpha\beta^{-1}\sup g$. 
In the latter case, blow-up of the delayed system alone may occur because the initial memory $v(0)$ is large.

As the analysis of this paper shows, the nature of finite-time blow-up in distributed delay differential equations depends strongly on the size of $v_{\max}$. 
When $v_{\max}=\infty$, finite-time singularities are limited to blow-up in which $u(t)$ diverges to infinity, as shown in Lemma \ref{lem:DDDEg-pro5}. 
On the other hand, when $v_{\max}<\infty$, the accumulated memory reaching the limit value ($v(t)\to v_{\max}$) can trigger blow-up. 
Assumption (V) introduced in this study captures the singular growth rate of the nonlinearity near $v\to v_{\max}$, and serves as a key to describing the diverse blow-up scenarios specific to delayed systems. 
A complete understanding of the threshold phenomena caused by the size of the initial memory when $v_{\max}<\infty$ is left as an open problem.

Moreover, the framework developed in this paper, based on subsolutions, supersolutions, and integral conditions, has made it possible to bring out the dynamical mechanisms specific to differential equation systems: the singularity caused by the finiteness of $v_{\max}$, and the separation of memory-induced blow-up through comparison with non-delayed equations. 
The main results and the framework of this paper are not limited to the distributed delay differential equation \eqref{eq:DDDEg-int1}; they may also be applicable to the case where the delay kernel is a more general gamma distribution, and to Volterra-type integral equations. 
We hope that this paper serves as a first step toward further generalization of these methods, the clarification of threshold phenomena caused by the size of the initial memory, and applications to more general delay and integral equations.

\section*{Acknowledgments}
This work was partially supported by JSPS KAKENHI Grant Number 25K17306. 
The author developed the overall framework and presented the main theorem and its proof using Claude (Anthropic). 
Claude (Anthropic) also assisted with the translation and language editing of the manuscript.
The author carefully reviewed, verified and substantially revised all materials used, and assumes full responsibility for the accuracy and final content of this paper.


\end{document}